\documentclass[11pt]{article}

\usepackage{graphicx}
\usepackage{amsmath}
\usepackage{amsfonts}
\usepackage{amssymb}
\usepackage[all]{xy}
\usepackage[french]{babel}

\newtheorem{theorem}{Th\'eor\`eme}
\newtheorem{lemma}{Lemme}
\newtheorem{proposition}{Proposition}

\newtheorem{definition}{D\'efinition}
\newtheorem{remark}{Remarque}
\newtheorem{example}{Exemple\/}

\newenvironment{proof}[1][D\'emonstration]{\textbf{#1.} }{\
  \rule{0.5em}{0.5em}}

\title{Contractions d'alg\`ebres de Jordan en dimension 2}

\author{J. M. Ancochea Berm\'udez\dag \footnote{e-mail:
    ancochea@mat.ucm.es}, R. Campoamor-Stursberg\dag \footnote{e-mail:
    rutwig@mat.ucm.es}\\ L. Garc\'{\i}a Vergnolle
  \dag \footnote{e-mail: lucigarcia@mat.ucm.es} \footnote{Le
    troisi\`eme auteur (L. G. V.) remercie la Fundaci\'on Ram\'on
    Areces qui finance sa bourse pr\'edoctorale.}, J. S\'anchez Hern\'andez\dag  \\
  \\
  \dag\ Dpto. Geometr\'{\i}a y Topolog\'{\i}a,\\ Facultad
  CC. Matem\'aticas U.C.M.\\Plaza de Ciencias 3, E-28040 Madrid\\
} \date{}

\begin{document}

\maketitle

\begin{abstract}
On d\'etermine les classes d'isomorphismes des alg\`ebres de
Jordan en dimension deux sur le corps des nombres r\'eels. En
utilisant des techniques d¡'Analyse Non Standard, on \'etudie les
propri\'et\'es de la vari\'et\'e des lois d'alg\`ebres de Jordan,
et aussi les contractions parmi ces alg\`ebres.
\end{abstract}
\medskip
Mots clefs: alg\`ebre de Jordan, rigidit\'e, contraction.
\section{D\'efinitions et propri\'et\'es pr\'eliminaires}
Le but de ce travail c'est d'\'etudier certaines propri\'et\'es de
la vari\'et\'e des lois d'alg\`ebres de Jordan en dimension 2.
D'abord, on classifie ces alg\`ebres sur le corps des nombres
r\'eels, et on introduit la notion de perturbation d'alg\`ebres de
Jordan en utilisant la th\'eorie des ensembles internes de Nelson
\cite{Ne}. Ceci, nous permet de d\'eterminer les composantes
ouvertes de la vari\'et\'e $J^2$, d'o\`u en r\'esulte que la
vari\'et\'e est form\'ee par trois composantes, deux de dimension
4 et une de dimension 2. Les autres alg\`ebres r\'esultent de
limites par contraction des alg\`ebres rigides d\'efinissant les
composantes ouvertes.
\begin{definition}
Une loi d'alg\`ebre de Jordan sur $\mathbb{R}$ est une application
bilin\'eaire sym\'etrique $\varphi: \mathbb{R}^{n}\times
\mathbb{R}^{n}\rightarrow \mathbb{R}^n$ qui v\'erifie l'identit\'e
\begin{equation}
\varphi\left(\varphi(X,X),\varphi(X,Y)\right)-\varphi\left(X,\varphi
\left(\varphi(X,X),Y\right)\right)=0,\quad
X,Y\in\mathbb{R}^{n}.\label{J1}
\end{equation}
On d\'esignera par $J^{n}$ l'ensemble des lois d'alg\`ebres de
Jordan sur $\mathbb{R}^{n}$.
\end{definition}

On appelle alg\`ebre de Jordan de dimension $n$ a toute paire
$(\mathbb{R}^{n},\varphi)$, o\`u $\varphi\in J^{n}$. De fa\c{c}on
naturelle, on peut aussi d\'efinir les notions d'id\'eal et
sous-alg\`ebre d'une alg\`ebre de Jordan \cite{Cr}.

\medskip

Rappelons qu'un sous-espace vectoriel $V$ est dit isotrope s'il
existe un vecteur non nul $v$ tel que $\varphi(v,w)=0$ pour tout
$w\in V$. Si le sous-espace $\mathbb{R}v$ est lui m\^{e}me
isotrope, on dit que $v$ est un vecteur isotrope.

\begin{definition}
Une alg\`ebre de Jordan sans \'el\'ements isotropes est dite
simple si elle ne poss\`ede pas d'id\'eaux  non nuls.
\end{definition}

\begin{proposition}
Toute sous-alg\`ebre d'une alg\`ebre de Jordan sans \'el\'ements
isotropes poss\`ede un \'el\'ement neutre.
\end{proposition}

Voir \cite{Ti} pour le d\'etail de la d\'emonstration.

\medskip

Soit $\left\{e_{1},\ldots,e_{n}\right\}$ une base de $\mathbb{R}^{n}$.
On peut identifier toute alg\`ebre de Jordan $\varphi$ avec ses
constantes de structure sur une base donn\'ee. De l'identit\'e
(\ref{J1}) on obtient que les coordon\'ees d\'efinies par
$\varphi(e_{i},e_{j})=a_{ij}^{k}e_k$ sont des solutions du syst\`eme:
\begin{equation}
  a_{ii}^{h}a_{kj}^{l}a_{lk}^{r} - a_{ii}^{h}a_{jk}^{l}a_{lh}^{r} -
  2a_{ij}^{h}a_{ik}^{l}a_{hl}^{r}  +
  2a_{il}^{r}a_{hj}^{l}a_{ik}^{h}=0,\quad 1\leq i,j,k,l,r,h\leq n.
  \label{J2}
\end{equation}

\section{Classification des alg\`ebres de Jordan r\'eelles en
dimension 2}

Dans le cas particulier de dimension 2, toute alg\`ebre de Jordan
est donn\'ee par les relations:
\begin{align}
\begin{split}
e_{1}\circ e_{1}&= a_{1}e_{1}+ a_{2}e_{2} \\
e_{2}\circ e_{2}&= b_{1}e_{1}+ b_{2}e_{2} \\
e_{1}\circ e_{2}&= c_{1}e_{1}+ c_{2}e_{2}.
\end{split}
\end{align}

Cette dimension a \'et\'e d\'ej\`a consider\'ee dans \cite{An}
d'un point de vue g\'eom\'etrique, et dans \cite{Sv} pour des
applications aux \'equations diff\'erentielles. On peut donc
exprimer la structure par une matrice de coefficients du type
\[\begin{pmatrix} a_{1} & a_{2}\\
b_{1} & b_{2}\\
c_{1} & c_{2}
\end{pmatrix}
\]

En d\'eveloppant le syst\`eme (\ref{J2}), on obtient les douze
\'equations suivantes:
\begin{align}
\begin{split}
a_{2}\left( c_{1}c_{2}\right) &= a_{2}\left(a_{2}b_{1}\right)\\
b_{1}\left( c_{1}c_{2}\right) &= b_{1}\left(a_{2}b_{1}\right)\\
a_{2}\left( a_{1}b_{1}+b_{2}c_{1}\right) &=
a_{2}\left(b_{1}c_{2}+c_{1}^{2}\right) \\
b_{1}\left(a_{1}b_{1}+b_{2}c_{1}\right) &= b_{1}\left(
b_{1}c_{2}+c_{1}^{2}\right)  \\
a_{2}\left(a_{2}c_{1}+c_{2}^{2}\right) &=
a_{2}\left(a_{1}c_{2}+a_{2}b_{2}\right)\\
b_{1}\left(a_{2}c_{1}+c_{2}^{2}\right) &=
b_{1}\left(a_{1}c_{2}+a_{2}b_{2}\right)\\
a_{1}c_{1}c_{2}+2a_{2}b_{1}c_{2} &=
2c_{1}c_{2}^{2}+a_{1}a_{2}b_{1} \\
a_{1}a_{2}c_{1}+3a_{1}c_{2}^{2}+2a_{2}b_{2}c_{2} &=
2a_{2}c_{1}c_{2}+2c_{2}^{3}+a_{1}^{2}c_{2}+a_{1}a_{2}b_{2}\\
2c_{1}^{2}c_{2}+2b_{1}c_{2}^{2}+a_{1}^{2}b_{1}+a_{1}b_{2}c_{1} &=
3a_{1}b_{1}c_{2}+2b_{2}c_{1}c_{2}+a_{1}c_{1}^{2}\\
2c_{1}c_{2}^{2}+2a_{2}c_{1}^{2}+a_{1}b_{2}c_{2}+a_{2}b_{2}^{2}&=
b_{2}c_{2}^{2}+2a_{1}c_{1}c_{2}+3a_{2}b_{2}c_{1} \\
2a_{1}b_{1}c_{1}+3b_{2}c_{1}^{2}+b_{1}b_{2}c_{2} &=
a_{1}b_{1}b_{2}+b_{2}^{2}c_{1}+2c_{1}^{3}+2b_{1}c_{1}c_{2}\\
2a_{2}b_{1}c_{1}+b_{2}c_{1}c_{2}&=
a_{2}b_{1}b_{2}+2c_{1}^{2}c_{2}.
\end{split}\label{SJ}
\end{align}

Il en r\'esulte que $J^{2}$ est une vari\'et\'e alg\'ebrique
immerse dans $\mathbb{R}^{6}$.

\begin{theorem}
Toute alg\`ebre de Jordan r\'eelle de dimension 2 est isomorphe
\`a une des alg\`ebres suivantes:
\begin{enumerate}
\item $\psi_{0}:\quad e_{1}\circ e_{1}=e_{1},\quad e_{2}\circ
e_{2}=e_{1},\quad e_{1}\circ e_{2}=e_{2}.$

\item $\psi_{1}:\quad e_{1}\circ e_{1}=e_{1},\quad e_{2}\circ
e_{2}=0,\quad e_{1}\circ e_{2}=e_{2}.$

\item $\psi_{2}:\quad e_{1}\circ e_{1}=0,\quad e_{2}\circ
e_{2}=e_{2},\quad e_{1}\circ e_{2}=0.$

\item $\psi_{3}:\quad e_{1}\circ e_{1}=e_{2},\quad e_{2}\circ
e_{2}=0,\quad e_{1}\circ e_{2}=0.$

\item $\psi_{4}:\quad e_{1}\circ e_{1}=e_{1},\quad e_{2}\circ
e_{2}=0,\quad e_{1}\circ e_{2}=\frac{1}{2}e_{2}.$

\item $\psi_{5}:\quad e_{1}\circ e_{1}=e_{1},\quad e_{2}\circ
e_{2}=-e_{1},\quad e_{1}\circ e_{2}=e_2.$
\end{enumerate}
De plus, l'alg\`ebre $\psi_{5}$ est une alg\`ebre de Jordan
simple.
\end{theorem}

Dans ce qui suit, le symbole $\psi_i$ d\'esignera les alg\`ebres
du th\'eor\`eme.

\medskip
Si $\varphi$ est une alg\`ebre de Jordan sans \'el\'ements
isotropes, d'apr\`es la proposition 1 il existe un \'el\'ement
unit\'e. Alors, on peut trouver une base de $\mathbb{R}^{2}$ telle
que la matrice de $\varphi$ soit de la forme
\begin{equation}
\begin{pmatrix}
1 & 0\\
a & b\\
0 & 1
\end{pmatrix}
 \label{M1}
\end{equation}
Soit $x=x_{1}e_{1}+x_{2}e_{2}$ un \'el\'ement g\'en\'erique. Alors
\[
\varphi(x,x)=(x_{1}^{2}+ax_{2}^{2})e_{1}+(2x_{1}x_{2}+bx_{2}^{2})e_{2}.
\]
Si l'alg\`ebre de Jordan $\varphi$ n'a pas d'\'el\'ements isotropes
alors l'\'equation pr\'ec\'edente \'equivaut \`a la condition
\begin{equation}
b^{2}+4a\neq 0.
\end{equation}

\begin{lemma}
Toute alg\`ebre de Jordan $\varphi\in J^{2}$ sans isotropie est
isomorphe \`a $\psi_{0}$ ou $\psi_{5}$.
\end{lemma}

\begin{proof}
Si $\varphi$ ne poss\`ede pas d'isotropie, alors on peut trouver
une base $\left\{e_{1},e_{2}\right\}$ telle que la loi s'exprime
par la matrice (\ref{M1}) avec $4a+b^{2}\neq 0$.
\begin{enumerate}

\item Si $4a+b^{2}>0$, on prend dans $\psi_{0}$ le changement
de base $e_{1}^{\prime}=e_{1},
e_{2}^{\prime}=\frac{b}{2}e_{1}+\sqrt{a+\frac{b^{2}}{4}}e_{2}$, et
sur cette base on obtient que
\[
e_{1}^{\prime}\circ e_{1}^{\prime}=e_{1}^{\prime},\quad
e_{2}^{\prime}\circ e_{2}^{\prime}=a e_{1}^{\prime}+b
e_{2}^{\prime},\quad e_{1}^{\prime}\circ
e_{2}^{\prime}=e_{2}^{\prime}.
\]

\item Si $4a+b^{2}<0$, on prend dans $\psi_{5}$ le changement
de base $e_{1}^{\prime}=e_{1},
e_{2}^{\prime}=\frac{b}{2}e_{1}+\sqrt{-(a+\frac{b^{2}}{4})}e_{2}$.
Les relations sur la base transform\'ee sont alors
\[
e_{1}^{\prime}\circ e_{1}^{\prime}=e_{1}^{\prime},\quad
e_{2}^{\prime}\circ e_{2}^{\prime}=a e_{1}^{\prime}+b
e_{2}^{\prime},\quad e_{1}^{\prime}\circ
e_{2}^{\prime}=e_{2}^{\prime}.
\]
\end{enumerate}
D'o\`u le r\'esultat.
\end{proof}

\begin{lemma}
Toute alg\`ebre de Jordan de dimension 2 avec isotropie et
\'el\'ement unit\'e est isomorphe \`a $\psi_{1}$.
\end{lemma}

Sans perte de g\'en\'eralit\'e on peut supposer que $e_{1}$ est
l'\'el\'ement unit\'e et $e_{2}$ l'\'el\'ement isotrope. Alors
\[
e_{1}\circ e_{1}=e_{1},\quad e_{2}\circ e_{2}=0,\quad e_{1}\circ
e_{2}=e_{2}.
\]

\begin{lemma}
Soit $\varphi\in J^{2}$ une alg\`ebre de Jordan avec isotropie
sans \'el\'ement unit\'e. Alors $\varphi$ est isomorphe \`a
$\psi_{2}, \psi_{3}$ ou $\psi_{4}$.
\end{lemma}

\begin{proof}
Si $\varphi$ est sans unit\'e et isotrope, on peut trouver une
base $\left\{e_{1},e_{2}\right\}$ telle que les relations soient
donn\'ees par
\[
e_{1}\circ e_{1}=a_{1}e_{1}+a_{2}e_{2},\quad e_{2}\circ
e_{2}=0,\quad e_{1}\circ e_{2}=c_{1}e_{1}+c_{2}e_{2}.
\]
Le syst\`eme d'\'equations (\ref{SJ}) est r\'eduit aux relations
\begin{equation*}
c_{1}=0,\; a_{2}c_{2}^{2}=a_{1}a_{2}c_{2},\;
3a_{1}c_{2}^{2}=2c_{2}^{3}+a_{1}^{2}c_{2}.
\end{equation*}
On obtient deux possibilit\'es non-\'equivalentes:

\begin{enumerate}

\item $\varphi=\begin{pmatrix}
a_{1} & a_{2}\\
0 & 0\\
0 & 0
\end{pmatrix}, a_{1},a_{2}\in \mathbb{R}$,

\item $\varphi=\begin{pmatrix}
a_{1} & 0\\
0 & 0\\
0 & \frac{a_{1}}{2}
\end{pmatrix}, a_{1}\in\mathbb{R}$,
\end{enumerate}

Dans le premier cas, si $a_{1}\neq 0$, le changement de base
$e_{1}^{\prime}=e_{2},\;
e_{2}^{\prime}=\frac{1}{a_{1}}\left(e_{1}+\frac{a_{2}}{a_{1}}e_{2}\right)$

montre que
\[
e_{1}^{\prime}\circ e_{1}^{\prime}=e_{1}^{\prime},\quad
e_{2}^{\prime}\circ e_{2}^{\prime}=a e_{1}^{\prime}+b
e_{2}^{\prime},\quad e_{1}^{\prime}\circ
e_{2}^{\prime}=e_{2}^{\prime},
\]
d'o\`u $\varphi\simeq \psi_{2}$. Si $a_{1}=0$, on consid\`ere
la changement $e_{1}^{\prime}=\frac{1}{\sqrt{|a_{2}|}}e_{1},\;
e_{2}^{\prime}=e_{2}$. \newline Pour le deuxi\`eme cas, il suffit
de consid\'erer le changement de base donn\'e par
$e_{1}^{\prime}=\frac{1}{a_{1}}e_{1},\; e_{2}^{\prime}=e_{2}$ pour
obtenir l'isomorphisme $\varphi\simeq \psi_{4}$.

Pour finir la preuve du th\'eor\`eme, il suffit de voir que
$\psi_{5}$ est simple. On peut facilement v\'erifier que
$\psi_{5}$ ne poss\`ede pas d'id\'eaux propres non triviaux.
\end{proof}

\begin{remark}
Il en r\'esulte que sur le corps des nombres complexes, les
alg\`ebres de Jordan $\psi_{5}$ et $\psi_{0}$ sont
isomorphes. On obtient alors un exemple d'une alg\`ebre simple
dont la complexifi\'ee n'est pas simple.
\end{remark}

\medskip

Si on consid\`ere l'action du groupe lin\'eaire $GL(2,\mathbb{R})$
sur la vari\'et\'e $J^{2}$, cette action est d\'efinie par:
\begin{equation*}
\left(\varphi,f\right) = f^{-1}\circ \varphi\circ
\left(f,f\right).
\end{equation*}

On notera par $\mathcal{O}(\varphi)$ l'orbite d'un \'el\'ement
$\varphi$ par cette action. Les \'el\'ements de l'orbite sont
\'evidemment les alg\`ebres de Jordan isomorphes \`a $\varphi$.

\section{Perturbations d'alg\`ebres de Jordan}

Dans ce paragraphe, on applique la th\'eorie des ensembles
internes (I.S.T) de Nelson \cite{Ne} \`a l'analyse de la
vari\'et\'e $J^{2}$.\footnote{On peut \'egalement consulter le
livre \cite{Lu} pour des applications pr\'ecises de l'analyse
non-standard aux alg\`ebres de Lie} Dans l'annexe on pr\'esente
les axiomes les plus importants utilis\'es dans ce travail. Dans
ce qui suit, on suppose que $n$ est standard, d'o\`u on en
d\'eduit que la vari\'et\'e $J^{n}$ est standard. Soient
$\varphi\in J^{n}$ un \'{e}l\'{e}ment standard et ${x,y}$
deux
vecteurs standard de $\mathbb{R}^{n}$. Alors $\varphi\left( {x}%
,{y}\right)  $ est standard. Si ces deux vecteurs sont
limit\'{e}s, alors $\varphi\left(  {x},{y}\right)  $
est aussi limit\'{e}.

\begin{definition}
Soit $\varphi$ une loi standard de $J^{n}$. Une perturbation
$\varphi_{0}$ de $\varphi$ est une loi d'alg\`{e}bre de Jordan sur
$\mathbb{R}^{n}$ infiniment proche \`{a} $\varphi$,
c'est-\`{a}-dire, pour tous les vecteurs
${x},{y\in}\mathbb{R}^{n},$ $\varphi\left(  {x}%
,{y}\right)  $ et $\varphi_{0}\left(
{x},{y}\right)  $ sont infinement proches. On notera
une perturbation par $\varphi\sim\varphi_{0}$.
\end{definition}

\begin{proposition}
Si $\varphi{\in} J^{n}$, alors l'ombre $\varphi_{0}$ de
$\varphi$ appartient \`{a} $J^{n}$.
\end{proposition}

\begin{proof}
Soient ${x},{y\in}\mathbb{R}^{n}$. Comme $\varphi$ est une loi
d'alg\`{e}bre de Jordan, on a
\[
\varphi\left(  \varphi\left(  x,x\right)  ,\varphi\left(
x,y\right)  \right) -\varphi\left(  x,\varphi\left(  \varphi\left(
x,x\right)  ,y\right) \right)  =0.
\]
Comme chaque partie est limit\'{e}e, en consid\'{e}rant l'ombre on obtient%
\[
\varphi_{0}\left(  ^{o}\varphi\left(  x,x\right)
,^{o}\varphi\left( x,y\right)  \right)  -\varphi_{0}\left(
x,^{o}\varphi\left(  \varphi\left( x,x\right)  ,y\right)  \right)
=0,
\]
alors
\[
\varphi_{0}\left(  \varphi_{0}\left(  x,x\right)
,\varphi_{0}\left(
x,y\right)  \right)  -\varphi_{0}\left(  x,\varphi_{0}\left(  ^{o}%
\varphi\left(  x,x\right)  ,y\right)  \right)  =0,
\]
et appliquant une fois de plus les propri\'et\'es de l'ombre en
agissant sur les applications lin\'eaires, il en r\'{e}sulte que
\[
\varphi_{0}\left(  \varphi_{0}\left(  x,x\right)
,\varphi_{0}\left( x,y\right)  \right)  -\varphi_{0}\left(
x,\varphi_{0}\left(  \varphi _{0}\left(  x,x\right)  ,y\right)
\right)  =0.
\]
\end{proof}

\begin{lemma}
Soit $V$ un espace vectoriel de $\mathbb{R}^{n}$ $\left(  n\text{
standard}\right)  $. Alors l'ombre $^{o}V$ est un sous-espace de
$\mathbb{R}^{n}$ de la m\^{e}me dimension.
\end{lemma}

\begin{proof}
Par d\'{e}finition, $^{o}V$ est le seul sous-ensemble standard de
$\mathbb{R}^{n}$ dont les \'{e}l\'{e}ments sont des ombres des
\'{e}l\'{e}ments limit\'{e}s de $V$. Comme l'addition et la
multiplication externe par des scalaires sont continues, $^{o}V$
est un sous-espace vectoriel. La m\'{e}thode d'orthogonalisation
de Graam-Schmidt nous donne une base orthogonale de $V$. M\^{e}me
si cette base n'est pas g\'{e}n\'{e}ralement standard, leurs
vecteurs sont limit\'{e}s et admettent une ombre de norme 1. Comme
ces ombres sont orthogonales, on en d\'{e}duit que $V$ et $^{o}V$
ont la m\^{e}me dimension.
\end{proof}

Comme cons\'{e}quence imm\'{e}diate de ce r\'{e}sultat on obtient
le r\'esultat suivant:

\begin{proposition}
Soit $\varphi_{0}$ une loi standard de $J^{n}$ et soit $\varphi$
une perturbation de $\varphi_0$. Alors l'ombre d'une
sous-alg\`{e}bre (resp. un id\'{e}al) de $\left(
\mathbb{R}^{n},\varphi\right)  $ est une sous-alg\`{e}bre (resp.
un id\'{e}al) de $\left( \mathbb{R}^{n},\varphi _{0}\right)  $ de
la m\^{e}me dimension.
\end{proposition}

\begin{lemma}
Dans les conditions pr\'{e}c\'edentes, si \ $\left(
\mathbb{R}^{n},\varphi _{0}\right)  $ est sans isotropie, alors la
perturbation $\left( \mathbb{R}^{n},\varphi\right)  $ est aussi
sans isotropie.
\end{lemma}

\begin{proof}
Comme $\left(  \mathbb{R}^{n},\varphi_{0}\right)  $ est une
alg\`{e}bre de Jordan standard sans isotropie, pour tout
$x\in\mathbb{R}^{n}-\{0\}$ standard,
$\varphi_{0}\left(  x,x\right)  $ est standard non nul et v\'{e}rifie%
\[
\varphi\left(  x,x\right)  \,\sim\,\varphi_{0}\left( x,x\right) .
\]
On en conclut que pour tout standard non nul la relation
$\varphi\left(  x,x\right)  \neq0$ est satisfaite. Par
continuit\'{e} de $\varphi$, la m\^{e}me propri\'{e}t\'{e} est
valable pour tout $x\in
\mathbb{R}^{n}$ qui ne soit pas infiniment petit. Pour $x\in\mathbb{R}%
^{n}-\{0\}$ on consid\`{e}re $\frac{x}{\left\|  x\right\|  }$, qui
est limit\'{e} non infiniment petit. Alors
\[
\varphi\left(  \frac{x}{\left\|  x\right\|  },\frac{x}{\left\|
x\right\|
}\right)  \neq0\text{, }%
\]
et par bilin\'eairit\'{e}
\[
\frac{\varphi\left(  x,x\right)  }{\left\|  x\right\|  ^{2}}\neq0,
\]
d'o\`{u} $\varphi\left(  x,x\right)  \neq0$.
\end{proof}

\begin{proposition}
  Soit $\varphi$ une perturbation de la loi d'alg\`ebre de Jordan
  standard $\varphi_{0}$ sur $\mathbb{R}^{n}$. Si
  $G_{0}=(\mathbb{R}^{n},\varphi_{0})$ est simple, alors
  $G=(\mathbb{R}^{n},\varphi)$ est simple.
\end{proposition}

\begin{proof}
D'apr\`es le lemme 5, $G$ est sans isotropie. Si $I$ est un
id\'eal propre de $G$, alors $^{o}I$ est un id\'eal de $G_{0}$ de
la m\^{e}me dimension. En particulier, si $G_{0}$ est simple, ceci
 implique que $^{o}I$ est r\'eduit \`a z\'ero, d'o\`u la
simplicit\'e de $G$.
\end{proof}

\begin{definition}
Une alg\`ebre de Jordan standard $\varphi_{0}$ est dite rigide si
toute perturbation est isomorphe \`a $\varphi_{0}$.
\end{definition}

Cette d\'efinition est la traduction au langage non-standard de la
notion classique de rigidit\'e. En effet, si toute perturbation de
$\varphi_{0}$ est isomorphe \`a $\varphi_{0}$, le halo de
$\varphi_{0}$ est contenu dans l'orbite
$\mathcal{O}(\varphi_{0})$. Ceci implique que l'orbite est
ouverte, et par le principe de transfert, on obtient
l'\'equivalence.

\begin{example}
Soit $\psi_{0}\in J^{2}$ la loi de Jordan d\'{e}finie dans le
th\'{e}or\`{e}me 1, il existe alors une base $\{e_{1},e_{2}\}$ sur
laquelle la loi est donn\'{e}e  par la matrice:
\[
\begin{pmatrix}
1 & 0\\
1 & 0\\
0 & 1\\
\end{pmatrix}
\]
Consid\'{e}rons une perturbation $\varphi \in J^{2}$ de
$\psi_{0}$, dans cette m\^{e}me base elle est donn\'{e}e par:
\[
\begin{pmatrix}
1+\epsilon_{1} & \epsilon_{2}\\
1+\epsilon_{3} & \epsilon_{4}\\
\epsilon_{5} & 1+\epsilon_{6}\\
\end{pmatrix},
\]
o\`{u} les $\epsilon_{1},\dots,\epsilon_{6}$ sont infiniments petits.\\
En faisant le changement de variables
$e_{1}^{\prime}=\frac{1}{1+\epsilon_{6}}e_{1},\quad
e_{2}^{\prime}=e_{2}$, on peut supposer que $\varphi$ est de la
forme:
\[ \begin{pmatrix}
1+\epsilon_{1} & \epsilon_{2}\\
1+\epsilon_{3} & \epsilon_{4}\\
\epsilon_{5} & 1\\
\end{pmatrix}.
\]

Si on impose les \'{e}quations qui d\'{e}finissent la
vari\'{e}t\'{e} $J^{2}$, on obtient la condition:
\begin{equation}
\epsilon_{5}=\epsilon_{2}(1+\epsilon_{3}). \label{ec1}
\end{equation}
Par ailleurs, comme $\psi_{0}$ est sans isotropie, $\varphi$ l'est
aussi et poss\`{e}de un \'{e}l\'{e}ment neutre de la forme $\alpha
e_{1}^{\prime}+\beta e_{2}^{\prime}$ qui v\'{e}rifie:
\begin{align*}
  \varphi(e_{1}^{\prime},\alpha e_{1}^{\prime}+\beta e_{2}^{\prime})
  &=e_{1}^{\prime},\\
  \varphi(e_{2}^{\prime},\alpha e_{1}^{\prime}+\beta
  e_{2}^{\prime})&=e_{2}^{\prime}.
\end{align*}
d'o\`u \begin{align}
\begin{split}
1&=\alpha (1+\epsilon_{1})+\beta \epsilon_{5},\\
0&=\alpha \epsilon_{2}+\beta,\\
0&=\alpha \epsilon_{5}+\beta (1+\epsilon_{3}),\\
1&=\alpha+\beta \epsilon_{4}.
\end{split}\label{ec2}
\end{align}
En rempla\c{c}ant~(\ref{ec1}) dans la troisi\`{e}me \'{e}quation
de~(\ref{ec2}) puis $\beta$ dans la quatri\`{e}me, le syst\`eme
est r\'eduit \`a:
\begin{align*}
0&=\beta+\alpha \epsilon_{2},\\
1&=\alpha (1-\epsilon_{2}\epsilon_{4}).\\
\end{align*}

On en d\'{e}duit pourtant que $\alpha \sim 1$ et que $\beta \sim
0$. Dans la base $\{\alpha e_{1}^{\prime}+\beta
e_{2}^{\prime},e_{2}^{\prime}\}$, on peut alors trouver $\epsilon$
et $\epsilon^{\prime}$ infiniment petits tels que $\varphi$ soit
de la forme:
\[
\begin{pmatrix}
  1 & 0\\
  1+\epsilon & \epsilon^{\prime}\\
  0 & 1\\
\end{pmatrix}.
\] Comme $(\epsilon^{\prime})^{2}+4(1+\epsilon)>0$, $\varphi$ est
isomorphe \`{a} $\psi_{0}$. On en conclut donc que $\psi_{0}$ est
rigide.

\end{example}

Les alg\`ebres de Jordan rigides sont donc d'un grand
int\'er\^{e}t pour l'\'etude de la vari\'et\'e des lois, parce que
les orbites sont des composantes ouvertes, et la cl\^{o}ture donne
des composantes alg\'ebriques de la vari\'et\'e.\footnote{Dans ce
sens, on trouve des
  analogies avec la th\'eorie des alg\'ebres de Lie~\cite{Go}.}

\medskip

Soit $\varphi_{0}\in J^{n}$ une loi standard. Soient Id $\in
GL\left( n,\mathbb{R}\right)  $ et $f\in\frak{gl}\left(
n,\mathbb{R}\right)  $ standard. Pour tout $\varepsilon$
infiniment petit, l'endomorphisme
$\mathrm{Id}$ $+\varepsilon f$ appartient au groupe $GL\left(  n,\mathbb{R}%
\right)  $, et alors%
\[
\left( \mathrm{Id}+\varepsilon f\right) ^{-1}\varphi_{0}\left( \left(
    \mathrm{Id}+\varepsilon f\right) ,\left( \mathrm{Id}+\varepsilon
    f\right) \right) =\varphi_{0}+\varepsilon\left(
  \delta_{\varphi_{0}}f\right) +\varepsilon^{2}\left( \Delta\left(
    \varphi_{0},f,\varepsilon\right) \right) ,
\]
o\`u
\[
\delta_{\varphi_{0}}f\left( x,y\right) :=\varphi_{0}\left( f\left(
    x\right) ,y\right) +\varphi_{0}\left( x,f\left( y\right) \right)
-f\left( \varphi_{0}\left( x,y\right) \right) ,\;\forall x,y\in
\mathbb{R}^{n}.
\]
Comme l'ombre de la droite qui joint $\varphi_{0}$ et un point
infiniment proche $\varphi_{0}^{\prime}$ est une droite tangente \`{a}
l'orbite $\mathcal{O}\left( \varphi_{0}\right) $ en $\varphi_{0}$,
l'espace tangent est donn\'{e} par
\[
T_{\varphi_{0}}\mathcal{O}\left(  \varphi_{0}\right)  =\left\{
\delta _{\varphi_{0}}f\;:\;f\in\frak{gl}\left( n,\mathbb{R}\right)
\right\}  .
\]

\begin{example}
  Soit $\varphi_{0}$ une loi standard sans isotropie. D'apr\`{e}s la
  classification, on peut trouver une base $\left\{
    e_{1},e_{2}\right\} $ sur laquelle la loi est donn\'{e}e par la
  matrice
\[
\psi_{0}=
\begin{pmatrix}
  1 & 0\\
  1 & 0\\
  0 & 1
\end{pmatrix}\text{ ou }
\psi_5=
\begin{pmatrix}
  1& 0\\
  -1 & 0\\
  0 & 1
\end{pmatrix}.
\]
En prenant un \'{e}l\'{e}ment $f=\begin{pmatrix}
  a & b\\
  c & d
\end{pmatrix}
\in\frak{gl}\left( 2,\mathbb{R}\right) $ on trouve que
\begin{align*}
  \delta_{\psi_{0}}f\left( e_{1},e_{1}\right)
  &=ae_{1}+be_{2},\\
  \delta _{\psi_{0}}f\left( e_{2},e_{2}\right) &=\left( 2d-a\right)
  e_{1}+\left( 2c-b\right)
  e_{2},\\
  \delta_{\psi_{0}}f\left( e_{1},e_{2}\right) &=be_{1}+ae_{2},
\end{align*}
et
\begin{align*}
  \delta_{\psi_5}f(e_1,e_1)&=ae_1+be_2,\\
  \delta_{\psi_5}f(e_2,e_2)&=(-2d+a)e_1+(2c+b)e_2,\\
  \delta_{\psi_5}f(e_1,e_2)&=-be_1+ae_2,
\end{align*}
et pourtant
\[
T_{\psi_{0}}f=\begin{pmatrix}
  a & b \\
  2d-a & 2c-b \\
  b & a
\end{pmatrix},\quad T_{\psi_5}f=\begin{pmatrix}
  a & b \\
  -2d+a & 2c+b \\
  -b & a \end{pmatrix}.
\]
La dimension des orbites sont en cons\'equence $\dim\mathcal{O}\left(
  \psi_{0}\right) =\dim\mathcal{O}\left( \psi_{5}\right)=4$. Comme ces
alg\`{e}bres sont rigides, on trouve que la vari\'{e}t\'{e} $J^{2}$
poss\`{e}de deux composantes ouvertes de dimension 4.
\end{example}

\subsection{Th\'eor\`eme de d\'ecomposition d'un point}

Dans ce paragraphe, on \'enonce un r\'esultat de M. Goze \cite{Go}
concernant la rigidit\'e des lois d'alg\`ebres non-associatives.
Ceci admet une application naturelle aux alg\`ebres de Jordan.

\begin{theorem}
Soit $M_{0}$ un point standard de $\mathbb{R}^{n}$ avec $n$
standard. Tout point $M$ infiniment proche de $M_{0}$ admet une
d\'ecomposition de la forme
\begin{equation*}
  M=M_{0}+\epsilon_{1}v_{1}+\epsilon_{1}\epsilon_{2}v_{2}+\ldots+
  \epsilon_{1}\dots\epsilon_{p}v_{p}
\end{equation*}
qui v\'erifie les conditions
\begin{enumerate}
\item $\epsilon_{1},\ldots,\epsilon_{p}$ sont des nombres r\'eels
  infiniment petits. \item $v_{1},\ldots,v_{p}$ sont des vecteurs
  standard de $\mathbb{R}^{n}$ lin\'eairement ind\'ependants.
\end{enumerate}
De plus, si
\begin{equation*}
  M=M_{0}+\epsilon_{1}^{\prime}v_{1}^{\prime}+\epsilon_{1}^{\prime}\epsilon_{2}^{\prime}
  v_{2}^{\prime}+\ldots+\epsilon_{1}^{\prime}\dots\epsilon_{q}^{\prime}v_{q}^{\prime}
\end{equation*}
est une autre d\'ecomposition de $M$ qui v\'erifie les conditions
pr\'ec\'edentes, alors $p=q$, et
\begin{enumerate}
\item $v_{i}^{\prime}=\sum_{j=1}^{i}a_{i}^{j}v_{j}$ dont les
  $a_{i}^{j}$ sont des nombres standard pour tout $i,j$ telles que
  $a_{i}^{i}\neq 0$.
\item
  $\epsilon_{1}\dots\epsilon_{i}=a_{i}^{i}\epsilon_{1}^{\prime}\dots\epsilon_{i}^{\prime}+
  a_{i+1}^{i}\epsilon_{1}^{\prime}\dots\epsilon_{i+1}^{\prime}+\ldots+a_{p}^{i}
  \epsilon_{1}^{\prime}\dots\epsilon_{p}^{\prime}$.
\end{enumerate}
Le nombre entier $p$ qui d\'ecrit la classe d'\'equivalence d'un
point $M$ est appell\'e longueur de $M$.
\end{theorem}
\begin{remark}
  Le th\'eor\`eme nous donne, en partant des points $(M,M_{0})$ de
  $\mathbb{R}^{n}$, avec $M_{0}$ standard infiniment proche \`a $M$,
  une r\'ef\'erence standard $(v_{1},\ldots,v_{n})$ telle que le
  drapeau $(V_{1},(V_{1},V_{2}),\ldots,(V_{1},\ldots,V_{n}))$ est
  compl\`etement d\'etermin\'e. Si on impose que $M$ et $M_{0}$ sont des
  points appartenant \`a une vari\'et\'e di\-ff\'e\-ren\-tia\-ble ou
  alg\'ebrique immerse dans $\mathbb{R}^{n}$, la r\'ef\'erence
  ant\'erieure s'adapte \`a la g\'eom\'etrie tangente en $M_{0}$. En
  particulier, $v_{1}$ est un vecteur du c\^{o}ne tangent \`a la
  vari\'et\'e. \cite{Go}
\end{remark}

\subsection{Rigidit\'e des lois d'alg\`ebres de Jordan}
Soit $\varphi$ une perturbation de la loi standard $\varphi_{0}\in
J^{n}$. Alors on peut d\'ecomposer $\varphi$ comme \cite{Go}:
\begin{equation}
  \varphi= \varphi_{0}+
  \epsilon_{1}\varphi_{1}+\epsilon_{1}\epsilon_{2}\varphi_{2}+\ldots+
  \epsilon_{1}\dots\epsilon_{p}\varphi_{p},\label{D1}
\end{equation}
o\`u $\epsilon_{i}$ sont des nombres infinit\'esimaux et
$\varphi_{1},\ldots,\varphi_{p}$ sont des applications bilin\'eaires
ind\'ependantes de $\mathbb{R}^{n}$. Il en r\'esulte que $\varphi\in
J^{n}$ si et seulement si
\begin{align}
\begin{split}
  \delta_{\varphi_{0}}\varphi_{1}&=0,\\
  \varepsilon_{2}\delta_{\varphi_{0}}\varphi_{2}+\varepsilon_{1}%
  \delta_{\varphi_{1}}\varphi_{0}+\varepsilon_{1}^{2}\varphi_{1}^{3}%
  +\varepsilon_{1}\varepsilon_{2}\delta_{\varphi_{1}}\varphi_{3}+\varepsilon
  _{1}^{2}\varepsilon_{2}\delta_{\varphi_{1}}\varphi_{2}+\qquad&\\
  +\varepsilon_{1}\varepsilon_{2}^{2}\delta_{\varphi_{2}}\varphi_{1}%
  +\varepsilon_{1}\varepsilon_{2}\varepsilon_{3}\left(
    \varphi_{0},\varphi _{1},\varphi_{2}\right)+\ldots+&\\
  +\varepsilon_{1}^{2}\varepsilon_{2}^{3}\varepsilon_{3}^{3}\ldots\varepsilon_{p-1}^{3}\varepsilon_{p}^{2}%
  \delta_{\varphi_{p}}\varphi_{p-1}+
  \varepsilon_{1}^{2}\varepsilon_{2}^{3}\varepsilon_{3}^{3}\ldots\varepsilon
  _{p}^{3}\varphi_{p}^{3}&=0,
\end{split}\label{K1}
\end{align}
o\`{u}
\begin{align}
  \varphi_{i}\circ\varphi_{j}\circ\varphi_{k}\left( X,Y\right)
  =&\varphi_{i}\left( \varphi_{j}\left( X,X\right)
    ,\varphi_{k}\left( X,Y\right)  \right)\nonumber\\
  &-\varphi_{i}\left( X,\varphi_{j}\left( \varphi
      _{k}\left(  X;X\right)  ,Y\right)  \right), \nonumber\\
  \varphi_{i}^{3} =&\varphi_{i}\circ\varphi_{i}\circ\varphi_{i},\nonumber\\
  \left( \varphi_{i},\varphi_{j},\varphi_{k}\right)=&\varphi_{i}%
  \circ\varphi_{j}\circ\varphi_{k}+\varphi_{j}\circ\varphi_{k}\circ\varphi
  _{i}+\varphi_{k}\circ\varphi_{i}\circ\varphi_{j}+\nonumber\\&+\varphi_{i}\circ\varphi
  _{k}\circ\varphi_{j}+
  \varphi_{j}\circ\varphi_{i}\circ\varphi_{k}+\varphi_{k}\circ\varphi_{j}%
  \circ\varphi_{i}\nonumber\\
  \delta_{\varphi_{i}}\varphi_{j} =&\varphi_{i}\left(
    \varphi_{i}\left( X,X\right) ,\varphi_{j}\left( X,Y\right) \right)
  +\varphi_{i}\left( \varphi_{j}\left( X,X\right)
    ,\varphi_{i}\left(  X,Y\right)  \right) +\nonumber\\
  &+\varphi_{j}\left( \varphi_{i}\left( X,X\right) ,\varphi_{i}\left(
      X,Y\right) \right) -\varphi_{j}\left( X,\varphi_{i}\left(
      \varphi_{i}\left( X,X\right) ,Y\right) \right)
  -\nonumber\\
  &-\varphi_{i}\left( X,\varphi_{j}\left( \varphi _{i}\left(
        X,X\right) ,Y\right) \right) -\varphi_{i}\left( X,\varphi
    _{i}\left( \varphi_{j}\left( X,X\right) ,Y\right)
  \right).\label{P1}
\end{align}

En effet, en substituant la d\'ecomposition (\ref{D1}) dans
(\ref{J1}), en divisant par $\epsilon_{1}$~\footnote{$\epsilon_{1}$
  \'etant non nul si $\varphi\neq \varphi_{0}$.} et en consid\'erant
l'ombre on obtient (\ref{K1}). On d\'esignera par $G(\varphi_{0})$
l'espace vectoriel
\begin{equation}
  G(\varphi_{0})=\left\{
    \varphi\in\mathbb{R}^{\frac{n^{2}+n^{3}}{2}}:\quad
    \delta_{\varphi_{0}}\varphi =0\right\}.
\end{equation}

\begin{lemma}
  Soient $\varphi_{1},\varphi_{2}\in J^{n}$. Le plan
  g\'{e}n\'{e}r\'{e} par $\left\{ 0,\varphi_{1},\varphi_{2}\right\} $
  est contenu dans $J^{n}$ si et seulement si
  $\delta_{\varphi_{1}}\varphi_{2}=\delta_{\varphi_{2}}\varphi_{1}=0.$
\end{lemma}

\begin{proof}
  Soit $\varphi=\lambda\varphi_{1}+\mu\varphi_{2},$
  $\lambda,\mu\in\mathbb{R}$.  Alors%
\[
\varphi\left( \varphi\left( x,x\right) ,\varphi\left( x,y\right)
\right) -\varphi\left( x,\varphi\left( \varphi\left( x,x\right)
    ,y\right) \right)
=\lambda^{2}\mu\delta_{\varphi_{1}}\varphi_{2}\left( x,y\right)
+\mu^{2}\lambda\delta_{\varphi_{2}}\varphi_{1}\left( x,y\right) ,
\]
d'o\`{u} l'affirmation.
\end{proof}

\begin{proposition}
  Si $\psi_{4}\in J^{2}$ est telle qu'il existe une base $\left\{
    e_{1},e_{2}\right\} $ de $\mathbb{R}^{2}$ dont la loi est
  donn\'{e}e par la matrice
\[
\psi_{4}=
\begin{pmatrix}
  1 & 0\\
  0 & 0\\
  0 & \frac{1}{2}%
\end{pmatrix},
\]
alors $\psi_{4}$ est rigide. De plus, $\dim \mathcal{O}(\psi_{4})=2$.
\end{proposition}

\begin{proof}
  Par transfert, il suffit de d\'{e}montrer la proposition pour
  $\varphi_{0}$ et $\left\{ e_{1},e_{2}\right\} $ standard. Soit
  $\varphi$ une perturbation de $\varphi_{0}$ et
\[
\varphi=\varphi_{0}+\varepsilon_{1}\varphi_{1}+\ldots+\varepsilon_{1}%
\dots\varepsilon_{6}\varphi_{6}%
\]
la d\'{e}composition correspondante. Soit
\[
G\left( \varphi_{0}\right) =\mathcal{L}\left\{ \begin{pmatrix}
    1 & 0\\
    0 & 0\\
    0 & \frac{1}{2}%
\end{pmatrix},
\begin{pmatrix}
  0 & 0\\
  0 & 1\\
  \frac{1}{2} & 0
\end{pmatrix}
\right\}  .
\]
Si $\varphi_{1}\in G\left( \varphi_{0}\right) $, alors on a $\delta
_{\varphi_{1}}\varphi_{0}=0$. D'apr\`{e}s le lemme, il en r\'{e}sulte
que $\varphi_{1}\in J^{2}$ et pourtant $\varphi_{1}^{3}=0$. Si on
divise la deuxi\`{e}me \'{e}quation de~(\ref{K1}) par
$\varepsilon_{2},$ en passant \`{a} l'ombre on obtient que
$\varphi_{2}\in G\left( \varphi _{0}\right) $. Le m\^{e}me
raisonnement nous montre que $\varphi_{2}%
^{3}=0,\;\delta_{\varphi_{1}}\varphi_{2}=\delta_{\varphi_{2}}\varphi
_{0}=\delta_{\varphi_{2}}\varphi_{1}=0$ et $\left( \varphi_{0},\varphi
  _{1},\varphi_{2}\right) =0$. En divisant maintenant par
$\varepsilon_{3}$, la m\'{e}thode montre que $\varphi_{3}\in G\left(
  \varphi_{0}\right) $. Pourtant, on obtient que toute perturbation
est au plus de longueur 2. Par ailleurs, dans la base $\left\{
  e_{1},e_{2}\right\} $ cette perturbation s'\'{e}crit comme%
\[
\begin{pmatrix}
1+\varepsilon & 0\\
0 & \varepsilon^{\prime}\\
\frac{\varepsilon^{\prime}}{2} & \frac{1+\varepsilon}{2}%
\end{pmatrix}
.
\]
Le changement de base donn\'{e} par
\[
\begin{pmatrix}
  1 & 0\\
  \frac{-\varepsilon^{\prime}}{1+\varepsilon} & 1+\varepsilon
\end{pmatrix}
\]
nous donne l'isomorphisme cherch\'{e}.
\end{proof}
\begin{theorem}
L'alg\`ebre de Jordan $\psi_{2}=
\begin{pmatrix}
  0 & 0\\
  0 & 1\\
  0 & 0
\end{pmatrix}
$ ne se perturbe pas sur $\psi_{5}=
\begin{pmatrix}
  1 & 0\\
  -1 & 0\\
  0 & 1
\end{pmatrix}
$. En particulier, $\psi_{5}$ ne se contracte pas sur $\psi_{2}$.
\end{theorem}
\begin{proof}
On consid\`ere la perturbation $\varphi=
\begin{pmatrix}
  \epsilon_{1} & \epsilon_{2}\\
  \epsilon_{3} & 1+\epsilon_{4}\\
  \epsilon_{5} & \epsilon_{6}
\end{pmatrix}
$ de $\psi_{2}$, o\`u $\epsilon_{i}\sim 0$ sont des
infinit\'esimaux. Le changement
$e_{2}^{\prime}=\frac{1}{1+\epsilon_{4}}e_{2}$ permet de supposer
que $\epsilon_{4}=0$. Si $\epsilon_{5}\neq 0$, l'\'equation
(\ref{SJ}) implique la relation
\[
\epsilon_{5}\left(2\epsilon_{1}\epsilon_{3}+3\epsilon_{5}-1-2\epsilon_{5}^{2}-2\epsilon_{3}\epsilon_{6}\right)=
\epsilon_{3}\left(\epsilon_{1}-\epsilon_{6}\right).
\]
Comme $\epsilon_{5}\neq 0$, cela \'equivaut \`a
\[
-1\sim
2\epsilon_{1}\epsilon_{3}+3\epsilon_{5}-1-2\epsilon_{5}^{2}-2\epsilon_{3}\epsilon_{6}=
\frac{\epsilon_{3}}{\epsilon_{5}}\left(\epsilon_{1}-\epsilon_{6}\right).
\]
Comme $\epsilon_{1}-\epsilon_{6}\sim 0$, le nombre
$\frac{\epsilon_{3}}{\epsilon_{5}}$ est infiniment grand, d'o\`u
$\frac{\epsilon_{5}}{\epsilon_{3}}\sim 0$. Le changement de base
\[
e_{1}^{\prime}=e_{1}-\frac{\epsilon_{5}}{\epsilon_{3}} e_{2}\sim
e_{1},\; e_{2}^{\prime}=\lambda e_{2}, \lambda\sim 1
\]
permet de supposer que $\epsilon_{5}=0$. Si $\epsilon_{3}\neq 0$,
le syst\`eme (\ref{SJ}) implique les relations $\epsilon_{2}=0$ et
$\epsilon_{1}=\epsilon_{6}$. Il faut distinguer deux cas:
\begin{enumerate}

\item Si $\epsilon_{1}\neq 0$, on consid\`ere la base donn\'ee par
\[
e_{1}^{\prime\prime}=\frac{1}{\epsilon_{1}}e_{1}^{\prime},\;
e_{1}^{\prime\prime}=e_{2}^{\prime}.
\]
Il en r\'esulte que
\[
e_{1}^{\prime\prime}\circ
e_{1}^{\prime\prime}=e_{1}^{\prime\prime},\;
e_{2}^{\prime\prime}\circ e_{2}^{\prime\prime}=
\epsilon_{1}\epsilon_{3} e_{1}^{\prime\prime}+
e_{2}^{\prime\prime},\; e_{2}^{\prime\prime}\circ
e_{2}^{\prime\prime}=e_{2}^{\prime\prime},
\]
et comme $\epsilon_{1}\epsilon_{3}+\frac{1}{4}>0$, la perturbation
est isomorphe \`a $\psi_{0}$ d'apr\`es le lemme 1.

\item Si $\epsilon_{1}=0$, la pertubation est isomorphe \`a
$\psi_{2}$. Ceci r\'esulte imm\'ediament de la classification.
\end{enumerate}

Si $\epsilon_{3}= 0$, alors la loi de la perturbation est donn\'ee
par $\varphi=
\begin{pmatrix}
\epsilon_{1} & \epsilon_{2}\\
0 & 1\\
0 & \epsilon_{6}
\end{pmatrix}
$. Dans ce cas, $I=(e_{2}^{\prime})$ est un id\'eal de
$\varphi$, et, par simplicit\'e, la perturbation n'est pas
isomorphe \`a $\psi_{5}$.
\end{proof}

\begin{proposition}
Les seules alg\`{e}bres rigides dans $J^{2}$ sont
$\psi_{0},\psi_{4}$ et $\psi_{5}$.
\end{proposition}

\begin{proof}
Soit $\varphi_{0}\in J^{2}$ une loi standard non isomorphe \`{a}
$\ $aucune
des lois%
\[
\begin{pmatrix}
1 & 0\\
1 & 0\\
0 & 1
\end{pmatrix}
,\ \begin{pmatrix}
1 & 0\\
0 & 0\\
0 & \frac{1}{2}%
\end{pmatrix}
,\ \begin{pmatrix}
1 & 0\\
-1&0\\
0& 1
\end{pmatrix}.
\]
Il existe une base standard $\left\{  e_{1},e_{2}\right\}  $ de $\mathbb{R}%
^{2}$ telle que $\varphi_{0}$ s'exprime d'une des formes suivantes%
\[
\text{a)\ }\varphi_{0}=
\begin{pmatrix}
1 & 0\\
0 & 0\\
0 & 1
\end{pmatrix}
,\;\text{b)\ }\varphi_{0}=
\begin{pmatrix}
0 & 0\\
0 & 1\\
0 & 0
\end{pmatrix}
 ,\;\text{c)\ }\varphi_{0}=
\begin{pmatrix}
0 & 1\\
0 & 0\\
0 & 0
\end{pmatrix}.
\]
Dans le cas a), on obtient la perturbation $\varphi$ donn\'{e}e
par
\[
\varphi=
\begin{pmatrix}
1 & 0\\
\varepsilon & 0\\
0 & 1
\end{pmatrix}
,
\]
qui n'est pas isomorphe \`{a} $\psi_{0}$ parce que $\varphi$
est sans isotropie. Si on est dans le cas b), alors on
consid\`{e}re la perturbation
\[
\varphi=
\begin{pmatrix}
\varepsilon & 0\\
0 & 1\\
0 & 0
\end{pmatrix}
 ,
\]
qui est aussi sans isotropie, et pourtant non isomorphe \`{a}
$\psi_{0}$. Finalement, dans le cas c) on peut consid\'{e}rer
la perturbation $\varphi=
\begin{pmatrix}
\varepsilon & 1\\
0 & 0\\
0 & 0
\end{pmatrix}
 $, qui n'est pas isomorphe \`{a} $\psi_{0}$ pour
\^{e}tre sans isotropie. La rigidit\'{e} de $\psi_{5}$ est
\'{e}vidente par simplicit\'{e}.
\end{proof}

\section{Contractions d'alg\`ebres de Jordan}

En analogie avec les alg\`ebres de Lie, on peut d\'efinir
formellement une notion de limite dans la vari\'et\'e des lois
$J^{n}$ de la fa\c{c}on suivante: Soit $\varphi\in J^{n}$ une loi
d'alg\`{e}bre de Jordan et soit $f_{t}\in GL\left(
n,\mathbb{R}\right)  $ une famille d'endomorphismes non-singuliers
d\'ependants du param\`etre continu $t$. Si pour tous $X,Y\in
J^{n}$ la limite
\begin{equation}
\varphi^{\prime}\left(  X,Y\right)  :=\lim_{t\rightarrow0}\,f_{t}^{-1}%
\circ\varphi\left(  f_{t}\left(  X\right)  ,f_{t}\left(  Y\right)
\right) \label{GW}
\end{equation}
existe, alors $\varphi^{\prime}$ v\'erifie les conditions
(\ref{J1}), et pourtant $\varphi^{\prime}$ est une loi d'alg\`{e}bre de Jordan,
appell\'{e}e contraction de $\varphi$ par $\left\{ f_{t}\right\}
$. En utilisant l'action du groupe $GL\left( n,\mathbb{R}\right) $
sur la vari\'{e}t\'{e} des lois d'alg\`{e}bres de Jordan, c'est
facile de voir qu'une contraction de $\varphi$ correspond \`{a} un
point de la cl\^{o}ture de l'orbite $\mathcal{O}\left(
\varphi\right) $. Il en r\'esulte en particulier que les
alg\`ebres rigides ne sont pas des contractions \cite{Ri}.

C'est \'{e}vident que le changement de base
$e_{i}^{\prime}=te_{i},\ i=1,..,n$  induit une contraction de
toute alg\`{e}bre de Jordan sur l'alg\`{e}bre de Jordan
ab\'elienne. De plus, toute contraction non-triviale $\varphi
\longrightarrow\varphi^{\prime}$ v\'{e}rifie l'in\'egalit\'{e} suivante%
\begin{equation}
\dim\mathcal{O}\left(  \varphi\right)  >\dim\mathcal{O}\left(
\varphi ^{\prime}\right)  .
\end{equation}
Ceci nous donne un premier crit\`{e}re pour classifier les
contractions des alg\`{e}bres de Jordan. On peut v\'{e}rifier sans
difficult\'{e} que pour les alg\`{e}bres du th\'{e}or\`{e}me~1,
les dimensions de l'orbite sont les
suivantes:%
\begin{align}
\begin{split}
\dim\mathcal{O}\left(  \psi_{5}\right)& =\dim\mathcal{O}\left(
\psi _{0}\right)  =4,\\ \dim\mathcal{O}\left( \psi_{1}\right)
&=\dim \mathcal{O}\left(  \psi_{2}\right)
=3,\\
\dim\mathcal{O}\left(  \psi _{3}\right) &=\dim\mathcal{O}\left(
\psi_{4}\right)  =2.
\end{split}
\end{align}

\begin{proposition}
Soient $\psi_{i}$ les alg\`{e}bres de Jordan du
th\'{e}or\`{e}me 1. Alors $\psi_{1},\psi_{2}$ et
$\psi_{3}$ sont les seules alg\`{e}bres de Jordan qui
apparaissent comme contraction d'une alg\`{e}bre de Jordan. Plus
pr\'{e}cisement:

\begin{enumerate}
\item $\psi_{1}$ est une contraction des alg\`{e}bres
$\psi_{0}$ et $\psi_{5},$

\item $\psi_{2}$ est une contraction de $\psi_{0},$

\item $\psi_{3}$ est une contraction de $\psi_{0},\psi_{1}%
,\psi_{2}$ et $\psi_{5}.$
\end{enumerate}
\end{proposition}

\begin{proof}
Comme les alg\`{e}bres $\psi_{0},\psi_{4}$ et $\psi_{5}$
sont rigides, elles ne sont pas une contraction d'une autre
alg\`{e}bre de Jordan.

\begin{enumerate}
\item  Soit dans $\psi_{5}$ la famille des transformations
lin\'{e}aires
\[
f_{t}\left(  e_{1}\right)  =e_{1},\;f_{t}\left(  e_{2}\right)
=te_{2}.
\]
Dans la base transform\'{e}e $\left\{  e_{1}^{\prime}=f_{t}\left(
e_{1}\right)  ,\;e_{2}^{\prime}=f_{t}\left(  e_{2}\right)
\right\}  $, la loi
de $\psi_{5}$ s'exprime par la matrice%
\[
\begin{pmatrix}
e_{1}^{\prime}\circ e_{1}^{\prime}\\
e_{2}^{\prime}\circ e_{2}^{\prime}\\
e_{1}^{\prime}\circ e_{2}^{\prime}%
\end{pmatrix}
 =
\begin{pmatrix}
1 & 0\\
-t^{2} & 0\\
0 & 1
\end{pmatrix}
.
\]
En appliquant la formule, on obtient que la limite pour
$t\rightarrow0$ de~(\ref{GW}) donne une alg\`{e}bre isomorphe
\`{a} $\psi_{1}$. De mani\`{e}re analogue, la contraction
$\psi_{0}\longrightarrow\psi_{1}$ est d\'{e}finie par les
transformations%
\[
f_{t}\left(  e_{1}\right)  =e_{1},\;f_{t}\left(  e_{2}\right)
=te_{2}.
\]

\item  En d\'{e}finissant les changements de base param\'etr\'es
\[
f_{t}\left(  e_{1}\right)  =te_{1},\;f_{t}\left(  e_{2}\right)  =\frac{1}%
{2}\left(  e_{1}+e_{2}\right)  ,
\]
on v\'{e}rifie facilement que la limite~(\ref{GW}) existe et donne
une contraction $\psi_{0}\longrightarrow\psi_{2}.$

\item  La contraction de $\psi_{5}$ sur $\psi_{3}$ est
donn\'{e}e par les transformations
\[
f_{t}\left(  e_{1}\right)  =\frac{\sqrt{t}}{2}e_{1}+\frac{\sqrt{t}}{2}%
e_{2},\;f_{t}\left(  e_{2}\right)  =te_{1}.
\]
La contraction $\psi_{1}\longrightarrow\psi_{3}$ est
d\'{e}finie par
\[
f_{t}\left(  e_{1}\right)  =t\left(  e_{1}+e_{2}\right)
,\;f_{t}\left( e_{2}\right)  =t^{2}e_{2}.
\]
La contraction $\psi_{0}\longrightarrow\psi_{3}$ est
donn\'{e}e par
\[
f_{t}\left(  e_{1}\right)  =t^{2}e_{2},\;f_{t}\left(  e_{2}\right)
=te_{1}.
\]
Finalement, les changements de base dans $\psi_{2}$
d\'{e}finis par
\[
f_{t}\left(  e_{1}\right)  =e_{1}+te_{2},\;f_{t}\left(
e_{2}\right)
=t^{2}e_{2}%
\]
induisent une contraction
$\psi_{2}\longrightarrow\psi_{0}$.
\end{enumerate}
L'alg\`ebre $\psi_{4}$, \'etant rigide, n'est pas une contraction
des autres alg\`ebres. Par ailleurs, cette alg\`ebre ne poss\`ede
que la contraction sur l'alg\`ebre ab\'elienne.
\end{proof}
Dans la figure~1, ou r\'esume les contractions des alg\`ebres de
Jordan obtenues, et aussi la contraction canonique sur l'alg\`ebre
ab\'elienne.

\begin{center}
\noindent
\[\xymatrix@C=1.2cm{ *{\begin{pmatrix} 1 & 0 \\ -1 & 0\\ 0 & 1
\end{pmatrix}}
\ar[rd] \ar@/_2pc/[rrdd] \ar@/_3pc/[rrddd] & &
*{\begin{pmatrix} 1 & 0\\ 1 & 0 \\ 0 & 1\end{pmatrix}}
\ar[ld] \ar[rd] \ar@/_1pc/[dd] \ar@/^4pc/[ddd] & & \\
& *{\begin{pmatrix} 1 & 0 \\ 0 & 0\\ 0 & 1
\end{pmatrix}} \ar[rd] \ar@/_1pc/[rdd] & &
*{\begin{pmatrix} 0 & 0 \\ 0 & 1\\ 0 & 0\end{pmatrix}}
\ar[ld] \ar@/^1pc/[ldd] &
 \\
 & & *{\begin{pmatrix} 0 & 1 \\ 0 & 0 \\ 0 & 0\end{pmatrix}} \ar[d] & &
 *{\begin{pmatrix} 1 & 0\\ 0 & 0\\ 0 &
     \frac{1}{2}\end{pmatrix}} \ar@/^1pc/[lld]\\
 & & *{\begin{pmatrix} 0 & 0\\ 0 & 0\\ 0 & 0\end{pmatrix}} & & }\]
{\sc Fig. 1} -- Contraction des alg\`ebres de Jordan de $J^{2}$.
\end{center}

\subsubsection*{Remerciements}

Les auteurs ont \'et\'e soutenus par le projet de recherche
MTM2006-09152 du Ministerio de Educaci\'on y Ciencia. Ils
t\'emoignent aussi leur gratitude \`a  A. Elduque, M. K. Kinyon et
P. Zusmanovich pour leurs commentaires et interventions
constructives lors de la r\'edaction de ce manuscrit.

\newpage
\section*{Annexe}
On adjoint au langage de Z.F.E (th\'{e}orie des ensembles de
Zermelo-Fraenkel avec l'axiome du choix) le symbole \lq st\rq{}
(lire standard), les formules qui d\'{e}coulent de la th\'{e}orie
de Z.F.E sont dites \lq\lq internes\rq\rq{} et les nouvelles
formules s'appelleront \lq\lq externes\rq\rq.  Les axiomes de
I.S.T. sont ceux de Z.F.E. restreints aux formules internes et les
axiomes de transfert, d'id\'{e}alisation et de standardisation que
l'on d\'{e}crit ci-dessous. On note que: $\forall^{st}x\,A(x)$
abr\`{e}ge la formule $\forall x\, [x \text{ standard }A(x)]$.
\begin{description}
\item{(T) --- }{\bf Axiome de Transfert}\\
  Soit $A(x,t_{1},\dots,t_{k})$ une formule interne sans autres
  variables libres que $x$ et les $t_{i}$, alors:
\[
\forall t_{1}\, , \dots,\, \forall t_{k}\, [\forall^{st} x\,
A(x,t_{1},\dots,t_{k})\, \leftrightarrow \, \forall x\,
A(x,t_{1},\dots,t_{k})].
\]
{\sl Cons\'{e}quences de l'{}axiome de Transfert:}
\begin{enumerate}
\item Pour d\'{e}montrer qu'une certaine propri\'{e}t\'{e} interne
(qui
  d\'{e}pend uniquement des variables standard) est vraie pour tout
  $x$, il suffit de la prouver pour tout $x$ standard.
\item Les ensembles d\'{e}finis de fa\c{c}on unique par une
formule
  interne qui d\'{e}pend de variables standard, sont standard. Les
  objets classiques $0, 1, 2, \dots, \pi, e, \varnothing, \mathbb{N},
  \mathbb{Q}, \mathbb{R}, \mathbb{C}, \mathbb{R}^{27},
  \mathbb{S}^{3},\dots$ sont donc standard ainsi que tous les objets
  construits \`{a} partir d'eux par une proc\'{e}dure interne.
\end{enumerate}
\item{(I) --- }{\bf Axiome d'id\'{e}alisation}\\
  Soit $B(x,y)$ une formule interne qui a au moins deux variables
  libres $x$ et $y$, alors:
\[
\forall^{st}z\text{ fini }\exists x\: \forall y\in z\: B(x,y) \:
\leftrightarrow \: \exists x\: \forall^{st} y \: B(x,y).
\]
{\sl Cons\'{e}quences de l'{}axiome d'id\'{e}alisation:}
\begin{enumerate}
\item Soit $\rho$ une relation binaire d\'{e}finie sur un ensemble
  standard $E$. Si pour tout ensemble standard fini $F\subset E$ il
  existe $v_{F}\in E$ reli\'{e} \`{a} tous les \'{e}l\'{e}ments $u\in
  F$, alors il existe $w\in E$ reli\'{e} \`{a} tous les
  \'{e}l\'{e}ments standard de $E$. On peut ainsi construire des
  objets \lq\lq id\'{e}aux\rq\rq{} par rapport \`{a} la math\'{e}matique
  classique.
\item En appliquant (I) \`{a} la formule $B(v,u)=(v\in
\mathbb{N}\,
  \wedge \, u\in \mathbb{N} \, \wedge \, v>u)$, on obtient la formule
  :
\[
\exists w\in \mathbb{N} \  \forall^{st} u\in \mathbb{N} \ w>u.
\]
On dit que cet \'{e}l\'{e}ment $w$ est infiniment grand,
\'{e}videmment $w$ n'est pas standard. De fa\c{c}on analogue,
l'axiome (I) fournit des r\'{e}els infiniment grands
($\forall^{st} y,\: y<\vert x \vert$) et des infiniment petits
($0$ et $\frac{1}{x}$ o\`{u} $x$ est infiniment grand).
\end{enumerate}

\begin{definition}
  $x\in \mathbb{R}$ est limit\'{e} si $\exists^{st} y \ \vert x
  \vert<y$. Deux r\'eels $x$ et $y$ sont infiniment proches si $\vert x-y \vert$
  est infiniment petit, on notera $x\sim y$.
\end{definition}

\item{(S) --- }{\bf Axiome de Standardisation}\\
Soit $C(z)$ une formule (interne ou externe) alors:
$$
\forall^{st} x \ \exists^{st} y\ \forall^{st} \ [z\in y \:
\leftrightarrow \: z\in x\: \wedge \: C(z)].
$$
{\sl Cons\'{e}quences de l'{}axiome de Standardisation:}
\begin{enumerate}
\item Soit $C(z)$ une formule et $X$ un ensemble standard,
l'axiome
  (S) nous dit qu'il existe un ensemble standard tel que tous ses
  \'{e}l\'{e}ments standard appartiennent \`{a} $x$ et v\'{e}rifient
  $C(z)$.
\item Pour tout r\'{e}el $x$ limit\'{e}, il existe un unique
r\'{e}el
  standard $ ^{0}x$, appel\'{e} ombre de $x$, tel que $x\sim \,
  ^{0}x$.\\
  Pour prouver l'existence de l'ombre d'un r\'{e}el
  limit\'{e} $x$, on applique (S) \`{a} l'ensemble des r\'{e}els et
  \`{a} la formule $t<x$, ceci nous donne un ensemble standard $E$ tel
  que: $\forall^{st} t\, (t\in E\,\leftrightarrow\,t \in
  \mathbb(R) \; \wedge \; t<x)$.  Comme $x$ est limit\'{e} et en
  utilisant (T), on d\'{e}duit que $E$ est une partie major\'{e}e non
  vide de l'ensemble des r\'{e}els et donc $E$ poss\`{e}de une
  borne sup\'{e}rieure $\sup E$, on
  v\'{e}rifie facilement que $\sup E \sim x$.\\
  De m\^{e}me, l'ombre $^{0}A$ de $A\subset \mathbb{R}$ est
  l'unique ensemble standard dont les \'{e}l\'{e}ments standard
  sont les ombres des \'{e}l\'{e}ments limit\'{e}s de $A$.\\
  Toute fonction $f:\mathbb{R} \rightarrow \mathbb{R}$ \`{a} valeurs
  limit\'{e}es sur les r\'{e}els standard admet une ombre $^{0}f$
  telle que:
\[
\forall^{st} x\in \mathbb{R} \quad (^{0}f)(x)=\,^{0}(f(x)).
\]
\end{enumerate}
\end{description}
Supposons $n$ standard, $\mathbb{R}^{n}$ et $\mathbb{C}^{n}$ sont
alors standard. On dira que $x\in \mathbb{R}^{n}$ est limit\'{e}
si chacune de ses composantes le sont et $x\in \mathbb{C}^{n}$ est
limit\'{e} si les parties r\'eelles et imaginaires de chacune de
ses composantes le sont. Ainsi, on peut g\'{e}n\'{e}raliser les
concepts d\'{e}finis ant\'{e}rieurement sur $\mathbb{R}^{n}$ et
sur
$\mathbb{C}^{n}$.\\
\par
En utilisant des ultrafiltres convenables, qui existent gr\^{a}ce
\`{a} l'axiome du choix, on prouve que I.S.T. est une extension
conservative de Z.F.E. Cela signifie que les propri\'{e}t\'{e}s
internes d\'{e}montr\'{e}es dans I.S.T. ont aussi une
d\'{e}monstration dans Z.F.E.


\begin{thebibliography}{99}

\bibitem{An} J. M. Ancochea Berm\'udez, Sobre la variedad de leyes
de \'algebras de Jordan, Publ. I.R.M.A. 216/P-124, 1983.

\bibitem{Go} M. Goze, Etude locale de la vari\'et\'e des lois
d'alg\`ebres de Lie, Th\`ese d'\'etat, Mulhouse, 1982.

\bibitem{Lu} R. Lutz, M. Goze, Non Standard Analysis. A practical
guide with applications. Lecture Notes in Mathematics 881.
Springer Verlag, 1981.

\bibitem{Cr} K. McCrimmon, A taste of Jordan algebras, Springer
Verlag, N.Y., 2004.

\bibitem{Ne} E. Nelson, Internal Set Theory. A new approach to
nonstandard analysis, Bull. Amer. Math. Society 83 (1977),
1165-1198

\bibitem{Ri} Nijenhuis A, Richardson R W, Cohomology and deformation
 of algebraic structures, Bull. Amer. Math. Society 70 (1964), 406-411.

\bibitem{Ti} A. Tillier, Quelques applications g\'eom\'etriques
des alg\`ebres de Jordan, Publ. Dep. Math. Lyon, 14, fasc. 3,
1977.

\bibitem{Sv} S. I. Svinolupov, I'ordanovy algebry i integriruemye
sistemy, Funkts. Analiz i ego Pril. 27, 40-53, 1993.

\end{thebibliography}
\end{document}